\documentclass[12pt]{amsart2000}

\usepackage{amsmath2000}
\usepackage{amsthm2000}
\usepackage{amsxtra2000}			 
\usepackage{amsfonts}
\usepackage{amssymb}
\usepackage{amscd2000}

\def\ZZ{{\mathbb Z}}
\def\NN{{\mathbb N}}

\def\RR{{\mathbb R}}
\def\QQ{{\mathbb Q}}
\def\PP{{\textbf P}}

\def\F{\mathcal{F}}

\def\O{\mathcal{O}}

\def\H{{H}}
\def\tensor{\otimes}
\def\iso{\simeq}

\DeclareMathOperator{\Pic}{Pic}

\DeclareMathOperator{\cl}{Cl}

\theoremstyle{plain}

\newtheorem{theorem}{Theorem}[section]
\newtheorem{proposition}[theorem]{Proposition}
\newtheorem{corollary}[theorem]{Corollary}
\newtheorem{lemma}[theorem]{Lemma}

\theoremstyle{definition}

\newtheorem{conjecture}[theorem]{Conjecture}
\newtheorem{conjecture/question}[theorem]{Conjecture/Question}

\theoremstyle{remark}

\begin{document}

\title{Vanishing theorems on toric varieties}

\author[M. Musta\c{t}\v{a}]{Mircea Musta\c{t}\v{a}}

\address{Department of Mathematics, University of California,
Berkeley, CA, 94720 and Institute of Mathematics of
the Romanian Academy}
\email{{\tt mustata@math.berkeley.edu}}

\thanks{ 
\newline \mbox{ } \mbox { }
}

\begin{abstract}
We use Cox's description for sheaves on toric varieties and results
about local cohomology with respect to monomial ideals to give 
a characteristic-free approach to vanishing results on toric varieties.
As an application, we give a proof of a strong version of Fujita's
Conjecture in the case of toric varieties. We also prove that every
sheaf on a toric variety corresponds to a module over the homogeneous
coordinate ring, generalizing Cox's result for the simplicial case.
\end{abstract}

\maketitle

%%----------------------------------------------------------
%%  INTRODUCTION    
%%----------------------------------------------------------
\section*{\bf Introduction}

   Our main goal in this article 
 is to give a characteristic free approach to vanishing
results on arbitrary toric varieties. We prove that the vanishing
of a certain cohomology group depending on a Weil divisor is implied
by the vanishing of the analogous cohomology group involving a higher
multiple of that divisor. When the variety is complete and the divisor
is $\QQ$-Cartier, one recovers in this setting a theorem due to
Kawamata and Viehweg. We apply these results to prove a strong form 
of Fujita's conjecture on a smooth complete toric variety.

Let $X$ be a toric variety and $D_1,\ldots,D_d$ the invariant Weil
divisors on $X$, so that $\omega_X\simeq\O_X(-D_1-\ldots-D_d)$. In the
first part of the paper we deduce the following  generalization
of the theorem of Kawamata and Viehweg, for toric varieties.

\begin{theorem}
 Let $D$ be an invariant Weil divisor on $X$ as above.
Suppose that we have $E=\sum_{j=1}^da_jD_j$, with $a_j\in\QQ$ and
$0\leq a_j\leq 1$ such that for some integer $m\geq 1$, $m(D+E)$
is integral and Cartier. If for some $i\geq 0$ we have 
 $H^i(\O_X(D+m(D+E)))=0$, then $H^i(\O_X(D))=0$. In particular,
if $X$ is complete and there is $E$ aforementioned such that $D+E$ is
$\QQ$-ample, then $H^i(\O_X(D))=0$ for every $i\geq 1$.
\end{theorem}

As a particular case of this theorem, we see that if for some 
$m\geq 1$ and $L\in\Pic(X)$ we have $H^i(L^m(-D_{j_1}-\ldots-
D_{j_r}))=0$, then $H^i(L(-D_{j_1}-\ldots-D_{j_r}))=0$. 
The cases $r=0$ and $r=d$ of this assertion were known to hold by reduction
to a field of positive characteristic. Over such a field $X$ is Frobenius
split and one concludes using arguments in Mehta and
Ramanathan in \cite{MR}. The fact that $X$ is Frobenius split will follow
also from our results.

Our method yields other vanishing results as well. 
For example, we prove that if $X$ is a smooth toric variety
and $L\in\Pic(X)$ is such that for some $m\geq 1$ and  $i\geq 0$ we have
$H^i(\Omega^j_X\otimes L^m)=0$, then $H^i(\Omega^j_X\tensor L)=0$.
By taking $X$ complete and $L$ ample, we thus recover a theorem of Bott,
Steenbrink and Danilov. 

In the second part of the paper we give some applications of vanishing 
theorems on toric varieties. Our main result is a proof
in the toric case of a strong version of the following
conjecture due to Fujita (see \cite{La}
for  discussion and related results).

\begin{conjecture}
 Let $X$ be a smooth projective variety of dimension 
$n$ and $L\in\Pic(X)$ an ample line bundle. Then $\omega_X\otimes L^{n+1}$
is globally generated and $\omega_X\otimes L^{n+2}$ is very ample.
\end{conjecture}

When the ample line bundle $L$ is generated by global sections,
an argument of Ein and Lazarsfeld \cite{EL} based on vanishing results
proves the conjecture over a field of characteristic zero. Under the
same hypothesis on $L$, the first assertion of the conjecture
is proved in arbitrary characteristic by Smith \cite{Sm}. 
On a smooth projective toric variety every ample line bundle is 
very ample
(see\cite{De}), so these  results prove the conjecture in this setting.

We give a direct proof of a strengthened form of the conjecture
in the case of a toric variety, where instead of
making a conclusion about a power $L^m$, we make a statement
about any line bundle $L$ satisfying
 conditions on the intersection numbers 
 with the invariant curves.
We are able to replace $\omega_X$ by the negative of the sum
of any set of $D_i$, and
also improve the bound by one in the case when
$X$ is not the projective space. More precisely, we prove:

\begin{theorem}
 Let $X$ be an $n$-dimensional complete smooth toric 
variety, $L\in\Pic(X)$ a line bundle and $D_1,\ldots,D_m$ distinct prime
invariant divisors.
\begin{enumerate}
\item If $(L\cdot C)\geq n$ for every invariant integral curve
$C\subset X$, then $L(-D_1-\ldots,-D_m)$ is globally generated, unless
$X\iso\PP^n$, $L\simeq\O(n)$ and $m=n+1$.
\item If $(L\cdot C)\geq n+1$ for every invariant integral curve
$C\subset X$, then $L(-D_1-\ldots-D_m)$ is very ample, unless
$X\iso\PP^n$, $L\simeq\O(n+1)$ and $m=n+1$.
\end{enumerate}
\end{theorem}

\medskip

To obtain these results we use Cox's notion of homogeneous coordinate
ring of $X$.
When the fan defining $X$ is nondegenerate (i.e., 
it is not contained in a hyperplane), this is
 a polynomial ring $S=k[Y_1,\ldots,Y_d]$
together with a reduced monomial ideal $B$ and with a
 grading in the class group of $X$ which is compatible with the $\ZZ^d$-grading
by monomials.
In general we need to slightly adjust this definition, but we leave
this generalization for the core of the paper.
 As in the case of projective space,
each graded $S$-module $P$ gives a quasi coherent sheaf $\widetilde{P}$
on $X$ and for each $i\geq 1$,
 the Zariski cohomology $H^i(X, \widetilde{P})$ is the degree zero part
of the local cohomology module $H^{i+1}_B(P)$. This idea
has been used in  \cite{EMS} 
to give an algorithm for the computation of cohomology of coherent sheaves
on a toric variety. 

Our basic result says that if $P$ is in fact $\ZZ^d$-graded and if the
multiplication by $Y_j$ is an isomorphism in certain $\ZZ^d$-degrees,
then the same is true for $H^i_B(P)$. The main example is $P=S$
in which we get that the multiplication
$$\nu_{Y_j} : H^i_B(S)_{\alpha}\longrightarrow H^i_B(S)_{\alpha+e_j}$$
is an isomorphism for every $\alpha=(\alpha_j)\in\ZZ^d$
 such that $\alpha_j\neq -1$.
In particular, $H^i_B(S)_{\alpha}$ depends only on the signs
of the components of $\alpha$.
This case was used in \cite{EMS} in order to
describe the support of $H^i_B(S)$. Similar results 
for the Ext modules appear also in \cite{Mu} and \cite{Ya}.
Our second example is that of the modules giving the sheaves
$\Omega^i_X$ on a smooth toric variety. Using this result and the relation 
between the local cohomology of a module and the Zariski cohomology
of the corresponding sheaf, we deduce the various vanishing theorems.

\bigskip

In the first section of the paper we summarize the construction in 
\cite{Co}
for the homogeneous coordinate ring suitably generalized to be 
applicable also to toric varieties defined by a degenerate fan.
All the results can be easily extended to this context. 
 We prove that every quasicoherent
sheaf on a toric variety comes from a graded module over the homogeneous
coordinate ring, generalizing the result in \cite{Co}
 for the simplicial case.
We describe the relation between the local cohomology of modules 
and the cohomology of the associated sheaves. This is used in the second 
section to prove the vanishing results described above.

In order to apply these results, we need numerical characterizations
for ampleness and numerical effectiveness for the toric case and in the
third section we provide these results.
In the simplicial case, a toric Nakai criterion is given in
\cite{Od}. We show that the result holds for
an arbitrary complete toric variety. 
We also prove that $L\in\Pic(X)$ is globally generated if and only if
$(L\cdot C)\geq 0$, for every integral invariant curve $C\subset X$.
In particular, we see that $L$ is globally generated if
and only if it is 
numerically effective. These results have been recently obtained also
by Mavlyutov in \cite{Ma}. We mention a generalization in a different
direction due to Di Rocco \cite{DR} who proved that on a smooth toric variety,
$L\in\Pic(X)$ is $k$-ample if and only if $(L\cdot C)\geq k$
for every invariant curve $C\subset X$.

 As a consequence of the above results, we deduce that $L$ is big and nef
if and only if there is a map $\phi : X\longrightarrow X'$ induced by a fan
refinement (therefore $\phi$ is proper and birational) and
$L'\in\Pic(X')$ ample such that $L\simeq\phi^*(L')$. This easily implies
the version of Kawamata-Viehweg vanishing theorem for nef and big
line bundles.

The fourth section is devoted to the above generalization (in this context)
of Fujita's Conjecture and some related results. The proof goes by induction 
on the dimension of $X$, by taking the restriction to the invariant prime
divisors. The result which allows the induction says that
for every $l\geq 1$, if 
$L$ is a line bundle such that $(L\cdot C)\geq l$
 for every invariant
integral curve $C\subset X$, then for every invariant prime
divisor $D$ and every $C\subset X$ aforementioned, $(L(-D)\cdot C)\geq l-1$.
From the case $l=1$ we see that if $L$ is ample, then $L(-D)$
is globally generated.
We conclude this section by proving a related result, which characterizes
the situation in which $L$ is ample and $D$ is a prime invariant divisor, 
but $L(-D)$ is not ample.

A well-known ampleness criterion (see, for example, \cite{Fu}) can
be interpreted as saying that on a complete toric variety $X$,
$L\in\Pic(X)$ is ample if
and only if it is globally generated and the map induced
by restrictions
$$H^0(L)\longrightarrow\bigoplus_{x\in X_0}H^0(L|_{\{x\}})$$
is an epimorphism, where $X_0$ is the set of fixed points of $X$.

In the last section we generalize this
property of ample line bundles under the assumption that
$X$ is smooth. We prove that the analogous
map is still an epimorphism if we replace $X_0$
 with any set of pairwise disjoint
invariant subvarieties. In this case, the blowing-up $\widetilde{X}$
of $X$ along the union of these subvarieties
 is still a toric variety and we obtain
 the required surjectivity by applying to $\widetilde{X}$ the 
results in the second section.

\medskip
\noindent
{\bf Acknowledgements.}
This work started from a joint project with David Eisenbud and Mike Stillman
to understand the cohomology of sheaves on toric varieties. It is a pleasure
to thank them for encouragement and generous  support. We are also
very grateful to
William Fulton, Robert Lazarsfeld and Sorin Popescu for useful discussions
and to Markus Perling for his comments on an earlier version of this paper.
Last but not least, we acknowledge the referee's numerous comments and
suggestions which greatly improved the quality of our presentation.

\section{The homogeneous coordinate ring of a toric variety}

Let $k$ be a fixed algebraically closed field of arbitrary characteristic.
We will use freely the definitions and results on toric varieties
from \cite{Fu}. We first review the notation we are going to use.

Let $N\simeq\ZZ^n$ be a lattice and $M={\rm Hom}(N, \ZZ)$ the dual lattice.
For a rational fan $\Delta$ in $V=N\otimes\RR$, we have
the associated toric variety $X=X(\Delta)$. For every $i\leq n$, the
set of cones in $\Delta$ of dimension $i$ is denoted by $\Delta_i$.
The torus $N\otimes_{\ZZ}k^*$ acts on $X$, and by an invariant subvariety
of $X$ we mean a subvariety which is invariant under this action.

The closed invariant subvarieties of $X$ of dimension $i$
are in bijection with the set $\Delta_{n-i}$.
For each cone $\tau\in\Delta$ we denote by $V(\tau)$ the corresponding
subvariety. Recall that $V(\tau)$ is again
a toric variety and $\tau_1\subset\tau_2$ if
and only if $V(\tau_2)\subset V(\tau_1)$. 
In particular, the prime invariant Weil divisors $D_1,\ldots,D_d$ on $X$
correspond to the one dimensional cones in $\Delta$.
If $X$ is smooth, then so is each $V(\tau)$.

  Let $V'$ be the vector space spanned by $\Delta$, $N'=N\cap V'$
and $M'={\rm Hom}(N',\ZZ)$ its dual lattice. We have an exact sequence:
$$0\longrightarrow M'\longrightarrow {\rm Div}_T(X)\longrightarrow {\rm Cl}
(X)\longrightarrow 0,$$
where ${\rm Div}_T(X)=\oplus_{i=1}^d\ZZ D_i\simeq\ZZ^d$ is the group
of invariant Weil divisors and ${\rm Cl}(X)$ is the class group of $X$.

We fix a decomposition $M\simeq M'\times\ZZ^e$, where $e$ is the codimension
of $V'$ in $V$. We correspondingly have a decomposition
$X\simeq X'\times (k^*)^e$, where $X'$ is the toric variety defined
by $\Delta$ in $N'$.

\bigskip

The homogeneous coordinate ring of $X$ was introduced by Cox
in \cite{Co}
in the case when the fan $\Delta$ is not degenerate, i.e., is not 
contained in a hyperplane. We slightly generalize this notion
in order to allow an arbitrary toric variety,
following the suggestion in \cite{Co}. We first
review some of the definitions and the results in that paper, all of
which can be easily generalized to this context.

For each $i$ with $1\leq i\leq d$ we introduce an indeterminate $Y_i$,
corresponding to the divisor $D_i$. We introduce also the 
indeterminates $Y_j$ with $d+1\leq j\leq d+e$, and the homogeneous
coordinate ring of $X$ is the ring $S=k[Y_1,\ldots,Y_d,Y_{d+1}^{\pm 1}
\ldots,Y_{d+e}^{\pm 1}]$. Note that the decomposition
$M\simeq M\times\ZZ^e$ corresponds to a decomposition
$k[M]\simeq k[M']\otimes k[Y_{d+1}^{\pm 1},\ldots,Y_{d+e}^{\pm 1}]$.
 
For every effective divisor $D=\sum_{i=1}^da_iD_i$, we write
$Y^D$ for the corresponding monomial $\prod_{i=1}^dY_i^{a_i}\in S$.   
 On the ring $S$ we have a fine grading,
the usual $\ZZ^{d+e}$-grading by monomials. However,
in this section we will consider exclusively a coarse ${\rm Cl}(X)$-grading
defined by 
$$\deg(\prod_{i=1}^{d+e}Y_i^{a_i})=[\sum_{i=1}^da_iD_i]\in{\rm Cl}
(X).$$

In the ring $S$ there is a distinguished ideal which is a reduced
monomial ideal. For each cone $\sigma\in\Delta$ we put
$D_{\hat {\sigma}}=\sum_{i;\tau_i\not\subset\sigma}D_i$,
the sum being taken
over the divisors corresponding to one dimensional cones
outside $\sigma$
and $Y^{\hat {\sigma}}=Y^{D_{\hat {\sigma}}}$.
If $\Delta_{\rm max}$ is the set of maximal cones in $\Delta$, 
then $B=(Y^{\hat {\sigma}}\,\vert\,\sigma\in\Delta_{\rm max})$.

As in the case of projective space, a graded $S$-module $P$
gives a quasicoherent sheaf on $X$ by the following procedure.
$X$ is covered by the affine toric varieties
$U_{\sigma}={\rm Spec}\,k[\sigma^{\vee}\cap M]$, for $\sigma\in\Delta$.
Using the above decomposition of $k[M]$ and the argument in \cite{Co},
we obtain canonical
isomorphisms $k[\sigma^{\vee}\cap M]\iso (S_{Y^{\hat {\sigma}}})_0$
for every $\sigma\in\Delta$, which are pairwise compatible.
Therefore if $P$ is a graded $S$-module, on the affine piece
$U_{\sigma}$ we can consider the quasicoherent sheaf
defined by $(P_{Y^{\hat {\sigma}}})_0$.
 These sheaves glue together to give
a quasicoherent sheaf $\widetilde{P}$ on $X$. In this way
we get an exact functor $P\longrightarrow\widetilde{P}$ from graded
$S$-modules to quasicoherent sheaves. If $P$ is finitely generated,
then $\widetilde{P}$ is coherent.

In particular, if $\alpha\in\cl(X)$, $\O(\alpha)$ is defined to
be $\widetilde {S(\alpha)}$. As in \cite{Co}, if 
$\alpha=[D]$, then there is a natural isomorphism 
$\O(\alpha)\iso\O(D)$. Moreover, we have an isomorphism
of graded rings
$$S\iso\bigoplus_{\alpha\in\cl(X)}H^0(X,\O(\alpha)).$$
For a quasicoherent sheaf $\F$, we put $\F(\alpha):=\F\otimes \O(\alpha)$.

{\bf Remark}. In general, if $P$ is a graded 
$S$-module, the natural morphism 
$\widetilde{P}\otimes 
\O(\alpha)\longrightarrow\widetilde {P(\alpha)}$ is not an
isomorphism. However, it is an isomorphism if $\alpha\in\Pic(X)$.
Indeed, by taking a graded free presentation of $P$, we can reduce ourselves
to the case when $P=S(\beta)$ for some $\beta=[E]$. Since
$\alpha=[D]$ with $D$ locally invertible, $\O(\alpha)$ is 
invertible and the fact that the morphism
$\O (D)\otimes\O (E)\longrightarrow\O (D+E)$ is an isomorphism follows
now directly from the definition.

\smallskip

We prove now that every quasicoherent sheaf is isomorphic
to $\widetilde{P}$ for some graded $S$-module $P$. This was proved in 
\cite{Co} under the assumption that
 $X$ is simplicial. With a slightly different definition
for the homogeneous coordinate ring it was proved
 more generally for toric varieties
with enough effective invariant divisors by Kajiwara in \cite{Ka}.

\begin{theorem}\label{sheaf}
For every toric variety $X$ and every quasicoherent
sheaf $\F$ on $X$, there is a graded $S$-module $P$ such that
$\F\iso\widetilde{P}$.
\end{theorem}

\begin{proof}
We take $P=\oplus_{\alpha\in\cl(X)}
H^0(X,\F(\alpha))$, which is clearly a graded $S$-module. For simplicity,
we will use the notation $P_{\sigma}=P_{Y^{\hat \sigma}}$.

For each $\sigma\in\Delta$, there are canonical maps
$$\phi_{\sigma}\,:\,(P_{\sigma})_0\longrightarrow H^0(U_{\sigma},\F),$$
defined as follows. If ${s/{Y^D}}\in(P_{\sigma})_0$ such that
$s\in H^0(X,\F(\alpha))$ and $D$ is an effective divisor with
$[D]=\alpha$ and ${\rm Supp}\,D\cap U_{\sigma}=\emptyset$, then
$1/{Y^D}$ defines a section in $H^0(U_{\sigma}, \O(-\alpha))$
and $\phi_{\sigma}({s/{Y^D}})=(1/Y^D)s$ is the image of
$(1/Y^D, s)$ by the canonical pairing
$$H^0(U_{\sigma}, \O(-\alpha))\times H^0(X, \F(\alpha))
\longrightarrow H^0(U_{\sigma}, \F).$$

These morphisms glue together to give 
$\phi\,:\,\widetilde{P}\longrightarrow\F$ (note that $\F$ is assumed
to be quasicoherent). We will prove that $\phi$ is an isomorphism
by proving that $\phi_{\sigma}$ is an isomorphism for each
$\sigma\in\Delta$.

\smallskip

We first show that $\phi_{\sigma}$ is a monomorphism. Suppose that
$\phi_{\sigma}({s/Y^D})=0$ for some 
$s\in H^0(X, \F(\alpha))$ and $D$ effective, $[D]=\alpha$.

We may assume that ${\rm Supp}\,D=\cup_{\tau_i\not\subset\sigma}
V(\tau_i)$, and in this case we will prove that there is an integer
$N\geq 1$ such that $Y^{ND}s=0$ in $H^0(X, \F(\alpha+N\alpha))$.
In fact, we will find for each $\tau\in\Delta$ an integer
$N_{\tau}$ such that $Y^{N_{\tau}D}s\vert_{U_{\tau}}=0$. Then it
is clear that $N=\sum_{\tau}N_{\tau}$ satisfies the requirement.

From now on, we fix also $\tau\in\Delta$. Since $\sigma\cap\tau$
is a face of $\tau$, we can write $\sigma\cap\tau=\tau\cap u^{\perp}$
for some $u\in\tau^{\vee}\cap M$. If for each $v\in M$,  
the corresponding element of $k[M]$
is denoted by $\chi^v$, we consider the principal divisor
$D_0={\rm div}(\chi^u)$. It is effective on $U_{\tau}$, where its
support coresponds to the one-dimensional cones $\tau_i\subset
\tau\cap\sigma$.

We consider the restrictions of all the sections from above to
$U_{\tau}$: $s\vert_{U_{\tau}}\in H^0(U_{\tau}, \F(\alpha))$,
$Y^D\vert_{U_{\tau}}\in H^0(U_{\tau}, \O(\alpha))$ and
$(1/Y^D)\vert_{U_{\sigma}\cap U_{\tau}}\in
H^0(U_{\sigma\cap\tau}, \O(-\alpha))$.

Since $\phi_{\sigma}({s/Y^D})=0$ in $H^0(U_{\sigma}, \F)$,
we have that $s\vert_{U_{\sigma}}=0\in H^0(U_{\sigma}, \F(\alpha))$,
as the image of $(Y^D, \phi_{\sigma}({s/Y^D}))$ by the canonical
pairing
$$H^0(U_{\sigma}, \O(\alpha))\tensor H^0(U_{\sigma}, \F)\longrightarrow
H^0(U_{\sigma}, \F(\alpha)).$$

In particular, we have $s\vert_{U_{\sigma\cap\tau}}=0$.
But $U_{\sigma}\cap U_{\tau}=U_{\sigma\cap\tau}\subset U_{\tau}$
is a principal affine subset defined by
$Y^{D_0}\in H^0(U_{\tau}, \O_X)$.
Therefore, we get an
integer $t\geq 1$ such that $Y^{tD_0}s=0$ in $H^0(U_{\tau}, \F(\alpha))$.

If $a_{\tau'}$ and $a^0_{\tau'}$ are the coefficients of 
$V(\tau')$ in $D$ and $D_0$, respectively, and $N_{\tau}$
is such that $N_{\tau}a_{\tau'}\geq ta^0_{\tau'}$ for every 
one-dimensional face $\tau'\subset\tau$ (by the form of $D$ and
$D_0$, we can choose such an $N_{\tau}$), then
$Y^{N_{\tau}D}s=0$ in $H^0(U_{\tau}, \F(\alpha+N_{\tau}\alpha))$.
This follows from the fact that if $\tau''$ is an one-dimensional cone with
 $\tau''\not\subset\tau$,
then $\O(V(\tau''))\vert_{U_{\tau}}$ is invertible
and $Y^{V(\tau'')}$ is an invertible section in it. This completes the
proof of the fact that $\phi_{\sigma}$ is a monomorphism.

\smallskip

We prove now that $\phi_{\sigma}$ is an epimorphism. Let $t\in
H^0(U_{\sigma}, \F)$, and let $D=\sum_{\tau_i\not\subset\sigma}
D_i$ and $\alpha=[D]$. 

Using an analogous argument, we see that for each $\tau\in
\Delta$, there is an integer $N_{\tau}$ such that
$Y^{N_{\tau}D}t\vert_{U_{\sigma\cap\tau}}\in
H^0(U_{\sigma\cap\tau}, \F(N_{\tau}\alpha))$ can be extended to a section
in $H^0(U_{\tau}, \F(N_{\tau}\alpha))$. Indeed, with the notation 
and arguments we used before, we first find $N'_{\tau}$
such that $Y^{N'_{\tau}D_0}t$ can be extended to $U_{\tau}$
and then find $N_{\tau}$, as claimed.

If we apply this to two cones $\tau_1$, $\tau_2\in\Delta$
and take $N\geq N_{\tau_1}$, $N_{\tau_2}$, we see that
$Y^{ND}t$ can be extended to both $U_{\tau_1}$ and
$U_{\tau_2}$, giving sections $t_1$ and $t_2$, respectively.
Since $(t_1-t_2)\vert_{U_{\sigma\cap\tau_1\cap\tau_2}}=0$,
by applying to $\tau_1\cap\tau_2$
 the argument we used to show that $\phi_{\sigma}$
is a monomorphism, we find
$N_{\tau_{12}}$ such that $Y^{N_{\tau_{12}}D}t_1
=Y^{N_{\tau_{12}}D}t_2$ on $U_{\tau_1}\cap U_{\tau_2}$.

This shows that for large enough $N$, we can extend
$Y^{ND}t\vert_{U_{\sigma\cap\tau}}$ to 
$t_{\tau}\in H^0(U_{\tau}, \F(N\alpha))$ for every $\tau\in\Delta$
such that $t_{\tau_1}\vert_{U_{\tau\cap\tau_2}}=
t_{\tau_2}\vert_{U_{\tau_1\cap\tau_2}}$
for every $\tau_1$, $\tau_2\in\Delta$. Therefore
$t$ is in the image of $\phi_{\sigma}$, which completes the proof. 
\end{proof}

Using the same argument as in \cite{Co}, we deduce the following corollary.

\begin{corollary}\label{coherent}
 For every toric variety $X$ and every coherent sheaf
$\F$ on $X$, there is a finitely generated $S$-module $P'$ such that
$\F\iso\widetilde {P'}$.
\end{corollary}

\begin{proof}
With the notation in the proof of Theorem~\ref{sheaf},
we have seen that
$$\phi_{\sigma}\,:\,(P_{\sigma})_0\longrightarrow H^0(U_{\sigma},\F)$$
is an isomorphism for every $\sigma\in\Delta$.

Since $\F$ is coherent, this implies that we can find a finitely generated
graded submodule $P'\subset P$ such that $(P'_{\sigma})_0
=(P_{\sigma})_0$ for every $\sigma\in\Delta$. It is clear that this
$P'$ satisfies the assertion of the corollary. 
\end{proof}

As in the case of projective space, the cohomology of the
sheaf $\widetilde{P}$ can be expressed as the local cohomology of 
the module $P$ at the irrelevant ideal $B$. 

\begin{proposition}\label{local_cohomology}
 Let $P$ be a graded $S$-module. Then there
exist an isomorphism of graded modules
$$H^{i+1}_B(P)
\iso\bigoplus_{\alpha\in\cl(X)}H^i(X, \widetilde {P(\alpha)}),$$
for every $i\geq 1$ and an exact sequence
$$0\longrightarrow H^0_B(P)\longrightarrow P\longrightarrow
\bigoplus_{\alpha\in\cl(X)}H^0(X, \widetilde {P(\alpha)})\longrightarrow
H^1_B(P)\longrightarrow 0.$$
\end{proposition}

\begin{proof}
$X$ is covered by the affine open subsets $U_{\sigma}$,
$\sigma\in\Delta_{\rm max}$, and we compute the cohomology
of $\widetilde{P}$ as Cech cohomology with respect to this cover.

On the other hand, we can compute the local cohomology of $P$
at $B$ using the direct limit of Koszul complexes on the powers
of the generators of $B=(Y_{\hat {\sigma}}\,\vert\,\sigma\in\Delta_{\rm max})$
(see \cite{Ei}, Appendix 4.1).

Since for $\sigma_1,\ldots ,\sigma_t\in\Delta_{\rm max}$,
$\bigcap_{i=1}^t U_{\sigma_i}=U_{\sigma}$, where $\sigma=
\bigcap_{i=1}^t\sigma_i$ and
$$H^0(U_{\sigma}, \widetilde {P(\alpha)})=(P(\alpha)_{Y^{\hat {\sigma}}})_0
=(P_{Y_{\hat {\sigma}_1},\ldots, Y_{\hat {\sigma}_t}})_{\alpha},$$
we conclude as in the case of the projective space
(see \cite{Ei}, Appendix 4.1).
\end{proof}

\smallskip

{\bf Note.} In the
 situation in Proposition~\ref{local_cohomology}, suppose that
$P$ is in fact a $\ZZ^{d+e}$-graded $S$-module, so that the corresponding
sheaf $\widetilde{P}$ is equivariant with respect to the torus action. 
In this case the local cohomology module $H^{i+1}_B(P)$ is  
$\ZZ^{d+e}$-graded, too, and under the isomorphism
in Proposition~\ref{local_cohomology} this finer decomposition of
$H^{i+1}_B(P)$ corresponds to the eigenspace decomposition of the
Zariski cohomology of the different twists $\widetilde{P(\alpha)}$.

\section{Vanishing theorems}

    We keep the notation from the previous section.
However, from now on we consider on $S$ the fine $\ZZ^{d+e}$ grading
by monomials and all $S$-modules are assumed to be $\ZZ^{d+e}$-graded.
Note that this implies that the associated sheaf
is equivariant with respect to the torus action.
The canonical basis of
 $\ZZ^{d+e}$ will
be denoted by $f_1,\ldots, f_{d+e}$.

For every subset $I\subset\ZZ$ and every graded $S$-module $P$,
we will say that $P$ is
 $I$-homogeneous if for every $\alpha=(\alpha_j)\in\ZZ^{d+e}$
with $\alpha_j\not\in I$, the multiplication by $Y_j$:
$$\nu_{Y_j}:P_{\alpha}\longrightarrow P_{\alpha+f_j}$$
is an isomorphism. Our main example is $S$, which is
obviously $\{-1\}$-homogeneous.

\begin{proposition}\label{homogeneous}
If $P$ is an $I$-homogeneous $S$-module,
then $H^i_B(P)$ is $I$-homogeneous.
\end{proposition}

\begin{proof}
We compute the local cohomology module
as the cohomology of a Cech-type complex (see, for example,
\cite{Ei}, Appendix 4.1). Let us temporarily denote the
generators of $B$ by $m_1,\ldots, m_t$. For a subset
$L\subset\{1,\ldots, t\}$, let $m_L$ be the least common multiple of
$\{m_l\,\vert\,l\in L\}$. Since it is enough to prove the assertion
at the level of complexes,
we have to check that for
every $\alpha\in\ZZ^{d+e}$ with $\alpha_j\not\in I$, the
multiplication by $Y_j$:
$$\mu_{Y_j}:(P_{m_L})_{\alpha}\longrightarrow (P_{m_L})_{\alpha+f_j}$$
is an isomorphism.

This is obvious if $Y_j\vert m_L$.
Suppose now that $Y_j\not| m_L$. Then the assertion is clear
once we notice that in this case, if $m/m^s_L\in P_{m_L}$, then
$\deg(m/m^s_L)_j=\deg(m)_j$, so that we can apply
the fact that $P$ in $I$-homogeneous. 
\end{proof}

We consider now an example of $\{-1,0\}$-homogeneous $S$-modules.
These are the modules which define the exterior
powers $\Omega^i_X$ of the cotangent sheaf. 
For simplicity, in this case we will assume that $X$ is smooth.

It is shown by Batyrev and Cox in \cite{BC} that
if the fan defining $X$ is nondegenerate, then the cotangent bundle
on $X$ appears in an Euler sequence:
$$0\longrightarrow\Omega^1_X\longrightarrow
\oplus_{j=1}^d\O _X(-D_j)\longrightarrow\O _X^{d-n}\longrightarrow 0.$$
In general, we have $X\iso X'\times (k^*)^e$ with $X'$ as above and
$\Omega^1_X\iso p_1^*(\Omega^1_{X'})\oplus\O_X^e$. Therefore we can
include $\Omega_{X'}^1$ in an exact sequence:
$$0\longrightarrow\Omega^1_X\longrightarrow (\oplus_{j=1}^d
\O_X(-D_j))\oplus\O_X^e\longrightarrow\O_X^{d-n+e}\longrightarrow 0.$$

We consider the graded morphism inducing the epimorphism
in the second exact sequence:
$$E=(\oplus_{j=1}^d S(-f_j))\oplus S^e\longrightarrow F=S^{d-n+e}.$$
For each $i\geq 1$, let
 $M_i$ be the kernel of the induced map $\wedge^iE\longrightarrow
\wedge^{i-1}E\tensor F$.

\begin{lemma}\label{omega}With the above notation, we have
 \item{\rm (i)} $\widetilde M_i\iso\Omega^i_X$.
\item{\rm (ii)}  $M_i$ is $\{-1,0\}$-homogeneous.
\end{lemma}

\begin{proof}
(i) The assertion 
follows easily from the above mentioned
result of Batyrev and Cox and the fact that in the Euler sequence all
the sheaves are locally free.

(ii) Since $M_i$ is a submodule of $\wedge^i E$, which is free,
the multiplication by $Y_j$ on $M_i$ is injective.

Let $\alpha=(\alpha_j)\in\ZZ^{d+e}$, $\alpha_j\not\in\{-1,0\}$. Since 
$\wedge^{i-1}E\tensor F$ is free, the surjectivity of the map
$\nu_{Y_j}\,:\,(M_i)_{\alpha}\longrightarrow (M_i)_{\alpha+f_j}$
follows from the surjectivity of the analogous map for $\wedge^iE$.
The latter is surjective since $\wedge^i E$ is a direct sum of
modules of the form $S(-f_{j_1}-\ldots-f_{j_r})$ with
$r\leq i$ and $j_1<\ldots<j_r$.
\end{proof}

\begin{proposition}\label{general_vanishing}
Let $X$ be an arbitrary toric variety.
\item{\rm (i)} If $P$ is a $\{-1,0\}$-homogeneous $S$-module
and $L\in\Pic(X)$ is such that $H^i(\widetilde{P}\tensor L^m)=0$
for some $i\geq 0$ and $m\geq 1$, then $H^i(\widetilde{P}\tensor L)=0$.
In particular, if $X$ is projective and $L\in\Pic(X)$ is ample, then
$H^i(\widetilde{P}\tensor L)=0$ for all $i\geq 1$. 
\item{\rm (ii)} Let $P$ be a $\{-1\}$-homogeneous $S$-module such that for 
every $\alpha\in\cl(X)$, $\widetilde {P(\alpha)}\iso\widetilde{P}\tensor
\O (\alpha)$. Suppose that 
$D\in {\rm Div}_T(X)$ and that there is
$E=\sum_{j=1}^d a_jD_j$ with $a_j\in\QQ$ and $0\leq a_j\leq 1$
such that $m(D+E)$ is integral and Cartier for some integer $m\geq 1$.
If $H^i(\widetilde{P}\tensor\O_X(D+m(D+E)))=0$
for some $i\geq 0$, then $H^i(\widetilde{P}\tensor\O_X(D))=0$.
In particular, if $X$ is projective and we have $E$ aforementioned such that
$D+E$ is $\QQ$-ample, then $H^i(\widetilde{P}\tensor\O_X(D))=0$
for all $i\geq 1$.
\end{proposition}

\begin{proof}
(i) If $L=\O(\alpha)$, 
$H^i(\widetilde{P}\tensor L)=H^i(\widetilde {P(\alpha)})$
(see the remark in the first section).
We will restrict ourselves to the case $i\geq 1$ in order
to apply the isomorphism in Proposition~\ref{local_cohomology}.
 When $i=0$,
one can give a similar argument using the exact sequence in that proposition.

As already mentioned, we have 
 $$H^i(\widetilde{P}\tensor L)\iso\bigoplus_{\underline{b}}
H^{i+1}_B(P)_{\underline{b}},$$
where the direct sum is taken over those $\underline{b}=(b_1,\ldots,b_{d+e})
\in\ZZ^{d+e}$ such that $[\sum_{i=1}^db_iD_i]=\alpha$.
Since by hypothesis $H^i(\widetilde{P}\tensor L^m)=0$, 
for every $\underline{b}$ with $[\sum_{i=1}^db_iD_i]=\alpha$
we have $H^i_B(P)_{m\underline{b}}=0$.
Proposition~\ref{homogeneous} implies that
$$H^i_B(P)_{\underline{b}}\simeq H^i_B(P)_{m\underline{b}},$$
which proves the first
assertion. 

In the case of an ample line bundle $L$
on a projective toric variety, we have
$H^i(\widetilde{P}\tensor L^m)=0$
for $i\geq 1$ and $m\gg 0$, so that we are in the previous situation.

\smallskip

(ii) We proceed similarly. Using our hypothesis on $P$
and Proposition~\ref{local_cohomology}, for every $i\geq 1$ we have
$$H^i(\widetilde{P}\tensor\O(\alpha))\simeq H^i(\widetilde{P(D)})
\iso\bigoplus_{\underline{b}}H^{i+1}_B(P)_{\underline{b}},$$
where the direct sum is taken over those 
$\underline{b}=(b_1,\ldots,b_{d+e})\in\ZZ^{d+e}$ 
such that $[\sum_{i=1}^db_iD_i]=[D]$.

Using again the hypothesis on $P$ and the fact that $m(D+E)$ is 
Cartier (see the remark in the first section), we get
$$H^i((P(D+m(D+E)))\widetilde{\,\,})=0.$$

We fix some $\underline{b}\in\ZZ^{d+e}$
 with $[\sum_{i=1}^db_iD_i]=[D]$. We have to prove that
$H^{i+1}_B(P)_{\underline{b}}=0$. If $\underline{b}'=
\underline{b}+m(\underline{b}+\underline{a})$, where
$\underline{a}=(a_1,\ldots,a_d,0\ldots,0)$, then
$[\sum_{i=1}^db'_iD_i]=[D+m(D+E)]$, and therefore
$H^{i+1}_B(P)_{\underline{b}'}=0$.

Proposition~\ref{homogeneous} implies that in order to complete the proof,
it is enough to show that $b_j\geq 0$ if
and only if $(m+1)b_j+ma_j\geq 0$.
This follows easily from the fact that $0\leq a_j\leq 1$. 
\end{proof}

We apply Proposition~\ref{general_vanishing}
 in conjunction with Lemma~\ref{omega} for
$P=M_i$ and for $P=S$.

\begin{theorem}\label{vanishing1}
\item{\rm (i)} {\rm (Bott-Steenbrink-Danilov)} If $X$ is a smooth
 toric variety and $L\in\Pic(X)$ is such that
$H^i(\Omega^j_X\tensor L^m)=0$ for some $m\geq 1$ and $i\geq 0$,
 then $H^i(\Omega^j_X\tensor L)=0$.
In particular, if $X$ is projective and $L\in\Pic(X)$ ample, then
$\H^i(\Omega^j_X\tensor L)=0$ for every $i\geq 1$.
\item{\rm (ii)}  Let $X$ be an arbitrary toric variety, $D\in {\rm Div}_T(X)$
and $E=\sum_{j=1}^da_jD_j$, with 
$a_j\in\QQ$ and $0\leq a_j\leq 1$ such that for some integer
 $m\geq 1$ we have $m(D+E)$ integral and Cartier. If  
$H^i(\O_X(D+m(D+E)))=0$, then
$H^i(\O_X(D))=0$. In particular, if $X$ is projective
and there is $E$ aforementioned such that $D+E$ is $\QQ$-ample,
then $H^i(\O_X(D))=0$ for all $i\geq 1$.
\end{theorem}

\smallskip

{\bf Remark.} As pointed out by the referee, in the case
$P=S$ the assertion in Proposition~\ref{homogeneous} can be proved
also via the combinatorial description of the cohomology
of a sheaf of fractional ideals (see for example \cite{KKMS}, pg. 42).
More precisely, the graded components $H_B^{i+1}(S)_{\alpha}$
and $H_B^{i+1}(S)_{\alpha+f_j}$ (or, equivalently, the corresponding
eigenspaces of $H^i(X,{\mathcal O}(\alpha))$ and
 $H^i(X,{\mathcal O}(\alpha+f_j))$)
can be described as simplicial cohomology groups of certain subsets of
$\RR^n$. The assertion can be proved by showing that these spaces
are homotopically equivalent.
Note that the case $P=S$ is enough to give
the statement of Proposition~\ref{vanishing1} (ii).

\smallskip

If $D=\sum_{j=1}^db_jD_j$ is a $\QQ$-divisor, we define
$$\lceil D\rceil :=\sum_{j=1}^d\lceil b_i\rceil D_j,$$
where for any real number $x$, $\lceil x\rceil$ is the integer
defined by $x\leq\lceil x\rceil < x+1$.
Similarily, we define
$$\lfloor D\rfloor:=\sum_{j=1}^d\lfloor b_j\rfloor D_j,$$
where for every $x$, $\lfloor x\rfloor$ is the integer defined by
$x-1<\lfloor x\rfloor\leq x$.
$K_X$ denotes the canonical divisor $-\sum_{j=1}^dD_j$
so that $\omega_X=\O(K_X)$.

\begin{corollary}\label{vanishing2}
Let $X$ be a projective toric variety.
\item{\rm (i)} {\rm (Kawamata-Viehweg)} If $D=\sum_{j=1}^db_jD_j$ is a 
$\QQ$-Cartier ample $\QQ$-divisor, then
$H^i(\O_X(K_X+\lceil D\rceil))=0$
for every $i\geq 1$.
\item{\rm (ii)}  If $D$ is as above, then 
$H^i(O_X(\lfloor D\rfloor))=0$
for every $i\geq 1$.
\item{\rm (iii)} Let $L\in\Pic(X)$ be an ample bundle. If $D_{j_1},
\ldots , D_{j_r}$ are distinct prime invariant divisors, then
$H^i(L(-D_{j_1}-\ldots-D_{j_r}))=0$
for every $i\geq 1$.
\end{corollary}

\begin{proof}
 All these are particular cases 
of Theorem~\ref{vanishing1} (ii).
\end{proof}

{\bf Remark.} In the proof of Fujita's Conjecture we will use
the assertion in Corollary~\ref{vanishing2} for smooth varieties.
As the referee pointed out, when $X$ is smooth
it is possible to prove this assertion
directly,
by induction on dimension and descending induction on $r$, as for
$r=d$ this is just Kodaira's vanishing theorem.

\smallskip

 As we mentioned in the Introduction, some particular cases of the above
 results can be proved by reducing the problem to a toric
variety $X$ over a field of
 positive characteristic $p$ and prove that such a variety
is Frobenius split. This means that if $F$ is the Frobenius morphism,
then the canonical morphism $\O_X\longrightarrow F_{*}\O_X$ has
a left inverse. With the description for the cohomology we used
above this can be seen as follows.

First of all, by embedding $X$ as an open subvariety of a
complete toric variety, we may suppose that $X$ is complete. Next,
by taking a toric resolution of singularities, we may suppose that $X$ 
is also smooth (see \cite{MR}). Moreover,
an argument in that paper shows that
in this case,
if $\dim(X)=n$, then $X$ is Frobenius split if and only if
the morphism
$$f\,:\,H^n(\omega_X)\longrightarrow H^n(\omega_X^p)$$
induced by the Frobenius morphism is not trivial.
But $H^n(\omega_X)\simeq H^{n+1}_B(S)_{(-1,\ldots,-1)}\iso k$, all
the other components being zero. On the other hand, 
$$H^n(\omega_X^p)\simeq\bigoplus_{[\sum(a_j+p) D_j]=0}\H^{n+1}_B(S)
_{\underline{a}}$$
 has by Proposition~\ref{homogeneous}
 the component $H^{n+1}_B(S)_{(-p,\ldots,-p)}$
canonically isomorphic with $\,\,\,\,\,\,\,\,$
 $H^{n+1}_B(S)_{(-1,\ldots,-1)}$ and therefore
with $k$. It is easy to see that via these identifications, the
corresponding component of $f$ is just the Frobenius morphism of $k$,
and therefore $f$ is nonzero.   

For a different approach to Frobenius splitting in the toric context
and other applications we refer to Buch, Thomsen, Lauritzen and Mehta
\cite{BTLM}.

\section{Ample and numerically effective line bundles}

Our main goal in this section is to give the condition for a line
bundle to be ample or nef (i.e., numerically effective) in terms of 
the intersection with the invariant curves. For ampleness, this is the
toric Nakai criterion which is proved in \cite{Od} for the smooth case
and is stated also for the simplicial case. We obtain also a similar
condition for the nef property, both
the results holding for arbitrary complete toric
varieties. In particular, we will see that on such a variety,
a line bundle is nef if
and only if it is globally generated. With a different
proof, these results have been obtained also by Mavlyutov in \cite{Ma}.
 We use the ideas 
in \cite{Od} together with the description for the intersection
with divisors in the non-smooth case from \cite{Fu}.

We will apply these results to show that a line bundle $L$ on $X$ 
which is big and nef is a pull-back of an ample line bundle on $X'$,
for a proper birational equivariant map of toric varieties 
$\phi:X\longrightarrow X'$. Recall that a line bundle $L$ on
$X$ is called nef if for every curve $C\subset X$, $(L\cdot C)\geq 0$.

\begin{theorem}\label{nef}
 If $X$ is a complete toric variety
and $L\in\Pic(X)$, the following are equivalent{\rm :}
\item{\rm (i)} $L$ is globally generated.
\item{\rm (ii)} $L$ is nef.
\item{\rm (iii)} For every invariant integral curve $C\subset X$,
$(L\cdot C)\geq 0$.  
\end{theorem}

\begin{proof}
 (i)$\Rightarrow$(ii) is true in general
and (ii)$\Rightarrow$(iii) follows from the definition.

We now prove the implication (iii)$\Rightarrow$(i). Let $D$ be an invariant
Cartier divisor such that $L\iso\O (D)$. Recall that there is a function
$\psi=\psi_D:N\tensor{\RR}\longrightarrow{\RR}$
associated with $D$ which is linear
on each cone $\sigma\in\Delta$. It is defined in the following way:
if $D\vert_{U_{\sigma}}={\rm div}(\chi^{-u_{\sigma}})\vert_{U_{\sigma}}$,
then $\psi\vert_{\sigma}=u_{\sigma}\vert_{\sigma}$ 
(the notation is that used in the first section).

A well-known result(see \cite{Fu}, Section 3.3)
says that $L$ is globally generated if
and only if $\psi$ is convex. 
Recall that $\dim\,(X)=n$.
To prove that $\psi$ is convex, it is enough to prove that for every
$\sigma_1$, $\sigma_2\in\Delta_n$ with $\dim(\sigma_1\cap
\sigma_2)=n-1$,
 $\psi\vert_{\sigma_1\cup\sigma_2}$ is convex, 
i.e., for every $x\in\sigma_1$, $y\in\sigma_2$ and
$t\in [0,1]$ such that $tx+(1-t)y\in\sigma_1\cup\sigma_2$,
we have $\psi(tx+(1-t)y)\geq t\psi(x)+(1-t)\psi(y)$.

It is clear, therefore, from the definition of $\psi$
 that it is enough to prove that for each
$\sigma_1$, $\sigma_2$ as above and each $D_i=V(\tau_i)$, with
$\tau_i\subset\sigma_2\setminus\sigma_1$ a one-dimensional cone,
$$u_{\sigma_2}(v_i)\leq u_{\sigma_1}(v_i),$$
where $v_i$ is the primitive vector of $\tau_i$.

Let $D=\sum_{j=1}^da_jD_j$. Note that by definition,
if $D_j=V(\tau_j)$, $\tau_j\subset\sigma$, then $u_{\sigma}(v_j)
=-a_j$.
For $\sigma_1$ and $\sigma_2$ as above, let $\tau=\sigma_1\cap\sigma_2$.
Our hypothesis gives
$(D\cdot V(\tau))\geq 0$.
By definition, $(D+{\rm div}(\chi^{u_{\sigma_1}}))\vert_{U_{\sigma_1}}=0$.
Therefore $$D+{\rm div}(\chi^{u_{\sigma_1}})
=\sum_{\tau_i\subset\sigma_2\setminus
\sigma_1}b_iD_i+\ldots ,$$
where we wrote down only the divisors corresponding to cones
in $\sigma_1\cup\sigma_2$. Since $a_i=-u_{\sigma_2(v_i)}$ for
$\tau_i\subset\sigma_2$, we get
$$b_i=u_{\sigma_1}(v_i)-u_{\sigma_2}(v_i),$$ 
if $\tau_i\subset\sigma_2\setminus\sigma_1$.

On the other hand, let us denote by $\overline{e}$ the
generator of the one-dimensional lattice $N/N_{\tau}$
such that the classes of the primitive vectors of $\tau_i$
for $\tau_i\subset\sigma_2\setminus\sigma_1$
 are positive multiples of $\overline{e}$. Here $N_{\tau}$
denotes the subgroup of $N$ generated by $N\cap\tau$.
If for every $\tau_i$ aforementioned we write $\overline{v}_i=c_i\overline{e}$,
then the intersection formula in \cite{Fu}, Section 5.1 shows that
$$(D+{\rm div}(\chi^{u_{\sigma_1}})\cdot V(\tau))=b_i/c_i,$$
for every $\tau_i\subset\sigma_1\setminus\sigma_2$.
Since $0\leq (D\cdot V(\tau))=b_i/c_i$ and $c_i>0$,
we deduce that
$b_i\geq 0$ for every $\tau_i$ aforementioned.
From the formula for $b_i$ we see that the proof
is complete.
\end{proof}

\smallskip

{\bf Remark.} The equivalence between (i) and (ii) above
can be deduced also from the result of Reid from \cite{Re}, which says that
every effective one dimensional cycle on $X$ is rationally equivalent
to an effective sum of invariant curves.

\begin{theorem}\label{ample}{\rm (Toric Nakai criterion)}
If $X$ is a complete toric variety, with $\dim(X) = n$ and
$L\in\Pic(X)$, then the following are equivalent{\rm :}
\item{\rm (i)} $L$ is ample.
\item{\rm (ii)} For every invariant integral curve $C\subset X$,
$(L\cdot C)>0$.
\end{theorem}

\begin{proof}
The proof of the relevant implication
(ii)$\Rightarrow$(i) is the same as the above proof for the implication
(iii)$\Rightarrow$(i). We have just to use the fact that $L=\O (D)$
is ample if and only if $\psi_D$ is strictly convex and to replace
all the inequalities by strict inequalities.
\end{proof}

\smallskip

Recall that a line bundle $L\in\Pic(X)$
 is called big if for a certain multiple $L^m$, the rational
map it defines: $\phi_{L^m}:X\longrightarrow\PP^N$ has the
image of maximal dimension $n=\dim\,(X)$.

\begin{proposition}\label{big}
\item{\rm (i)} If $X$ is a complete toric variety 
of dimension $n$ and $L\in\Pic(X)$   
is a line bundle which is globally generated and big, then
$\dim\,\phi_L(X)=n$.
\item{\rm (ii)} $L\in\Pic(X)$ is globally generated and big if and only if
there is a fan $\Delta'$ such that $\Delta$ is a refinement of
$\Delta'$ and $L'\in\Pic(X')$ ample, where $X'=X(\Delta')$,
and that if $f:X\longrightarrow X'$ is the map induced by the refinement,
$f^*(L')\iso L$.
\end{proposition}

\begin{proof}
 Let us fix an invariant Cartier divisor $D$
such that $L\iso\O (D)$. If $\psi_D$ is the function which appeared
in the proof of Theorem~\ref{nef}, it defines an associated convex polytope
$$P_D=\{u\in M\tensor\RR\,\vert\,
 u\geq\psi_D\,{\rm on}\,N\tensor\RR\}.$$
If $L$ is globally generated, then $\dim\,\phi_L(X)=\dim\,P_D$
(see \cite{Fu}, Section 3.4). But 
$P_{mD}=mP_D$, so that $\dim\,\phi_L(X)=\dim\,\phi_{L^m}(X)$, which
completes the proof of (i).

Since a map as in (ii) is birational, the ``if'' part in (ii) is trivial.
Let us suppose now that $L$ is globally generated and big. By the
above argument, $P=P_D$ is an $n$-dimensional convex polytope.
Such a polytope defines a complete fan $\Delta'$ and an ample
Cartier divisor $D'$ on $X'=X(\Delta')$. The cones in $\Delta'$
are in a one-to-one correspondence, reversing inclusions,
 with the faces of $P$:
for a face $Q$ of $P$ we have the cone 
$$C_Q=\{v\in N\tensor\RR\,\vert\,
\langle u,v\rangle\,\leq\,\langle u',v\rangle
\,{\rm for}\,{\rm all}\,u\in Q,u'\in P\}.$$
For every $\sigma\in\Delta_n$, $u_{\sigma}$ is a vertex of $P$.
Indeed, it is the intersection of $P$ with 
$$\{u\in M\tensor\RR\mid \langle u,v_i\rangle\geq\psi_D(v_i)\,
{\rm for}\,v_i\in\sigma\}.$$
In fact, every vertex of $P$ is of this form. 
Indeed, if $u_0$ is a vertex of $P$, then there
is $v\in N\tensor\RR$ such that $\langle u_0,v\rangle
<\langle u,v\rangle$, for all $u\in P\setminus\{u_0\}$.
In particular, we have $\psi_D(v)=\langle u_0,v\rangle$.
If $\sigma\in\Delta_n$ is such that $v\in\sigma$, then
$\langle u_0,v\rangle=\langle u_{\sigma},v\rangle$,
so that $u_0=u_{\sigma}$.

Now it is easy to check that
 $$C_{u_{\sigma}}=\bigcup_{\tau\in\Delta_n, u_{\tau}=u_{\sigma}}\tau.$$
Therefore $\Delta$ is a refinement of $\Delta'$.
Moreover, the ample divisor $D'$ on $X'$ is defined by
$$\psi_{D'}(v)={\rm min}_{\sigma\in\Delta_n}
\langle u_{\sigma}, v\rangle
=\psi_D(v).$$
It follows that if $f:X\longrightarrow X'$ is the map induced by the
refinement, $f^*(D')= D$, which completes the proof. 
\end{proof}

It is easy to see that using the results of this section, we can 
extend the form of the Kawamata-Viehweg theorem we obtained
 in the previous section to the case of a divisor which is big and nef. 
For the proof, however, we have to assume that 
the divisor is Cartier.

\begin{theorem}{\rm (Kawamata-Viehweg)}\label{kawamata}
If $X$ is a projective toric
variety and $L\in\Pic(X)$ is a line bundle which is nef and big,
then $\H^i(\omega_X\tensor L)=0$  for every $i\geq 1$.
\end{theorem}

\begin{proof}
Since $L$ is a line bundle, the duality
theorem gives
$$\H^i(X,\omega_X\tensor L)\iso\H^{n-i}(X,L^{-1}),$$
where $n=\dim(X)$ (see \cite{Fu}, Section 4.4).

Using Theorem~\ref{nef} and Proposition~\ref{big}, we get a morphism
$f\,:\,X\longrightarrow X'$, induced by a fan refinement, and 
$L'\in\Pic(X')$ ample such that $f^*(L')\iso L$. But then
$$\H^{n-i}(X, L^{-1})\iso\H^{n-i}(X', L'^{-1})
\iso\H^i(X',\omega_{X'}\tensor L')=0,$$
by Corollary~\ref{vanishing2}. 
\end{proof}

\begin{corollary}\label{base_locus}
 Let $X$ be a complete toric variety and
$L$ a line bundle on $X$. If  the base locus of $L$
is nonempty, then it contains an integral invariant curve $C\subset X$.
\end{corollary}

\begin{proof}
This is an immediate consequence of 
Theorem~\ref{nef}, since for an integral curve $C\subset X$,
if $C$ is not contained in the base locus of $L$,
then $(L\cdot C)\geq 0$.
\end{proof}

\section{Fujita's conjecture on toric varieties}

The main result of this section is the following strong form of
Fujita's Conjecture in the toric case.

\begin{theorem}\label{fujita}
 Let $X$ be an $n$-dimensional projective smooth
toric variety, $L\in\Pic(X)$ a line bundle and $D_1,\ldots,D_m$
distinct prime invariant divisors.
\item{\rm (i)} If $(L\cdot C)\geq n$ for every invariant intergral curve
 $C\subset X$,
then $L(-D_1-\ldots-D_m)$ is globally generated, unless
$X\iso\PP^n$, $L\iso\O(n)$ and $m=n+1$.
\item{\rm (ii)} If $(L\cdot C)\geq n+1$ for every
invariant integral curve $C\subset X$,
then $L(-D_1-\ldots-D_m)$ is very ample, unless $X\iso\PP^n$,
$L\iso\O(n+1)$ and $m=n+1$.
\end{theorem}

In particular, we have the following corollary.

\begin{corollary}\label{classic_fujita}
Let $X$ be an $n$-dimensional projective smooth
toric variety and $L\in\Pic(X)$.
\item{\rm (i)} If $(L\cdot C)\geq n$ for every
invariant integral curve $C\subset X$,
then $\omega_X\tensor L$ is globally generated, unless
$(X, L)\iso (\PP^n, \O(n))$.
\item{\rm (ii)} If $(L\cdot C)\geq n+1$ for every
invariant integral curve $C\subset X$,
then $\omega_X\tensor L$ is very ample, unless $(X, L)\iso (\PP^n, \O(n+1))$.
\end{corollary}

We prove Theorem~\ref{fujita}, using the numerical conditions for $L$
to be globally generated or ample, as well as the vanishing result  
in Corollary~\ref{vanishing2} (iii).
 The proof goes by induction on the dimension of $X$,
based on the following proposition.

\begin{proposition}\label{induction}
 Let $X$ be a projective smooth toric variety
with $\dim(X)=n$, $L\in\Pic(X)$ and $l\geq 1$ an integer. If
$(L\cdot C)\geq l$ for every invariant integral curve
$C\subset X$, then for every
prime invariant divisor $D$ and every $C$ aforementioned,
 $(L(-D)\cdot C)\geq l-1$.
\end{proposition}

We first deal with the case $l=1$ of this proposition in the lemma below.

\begin{lemma}\label{l1}
Let $X$ be a projective smooth toric variety,
$\dim(X)=n$. If $L\in\Pic(X)$ is ample and $D$ is an invariant prime
divisor, then $L(-D)$ is globally generated.
\end{lemma}

\begin{proof}[Proof of Lemma~\ref{l1}]
 We prove the lemma by induction on $n$.
For $n=1$, $X=\PP^1$ and the assertion is clear. 
If $n\geq 2$ and $L(-D)$ is not globally generated, since the
base locus of $L(-D)$ is invariant, we can choose a fixed point
$x$ in this locus.

Let $D'$ be a prime divisor distinct from $D$ and containing $x$.
By Corollary~\ref{vanishing2} (iii) , the restriction map
$$\H^0(L(-D))\longrightarrow\H^0(L(-D)\vert_{D'})$$
is surjective. On the other hand, $D'$ is a smooth toric variety of 
dimension $n-1$ and $D\cap D'$ is either empty or
a prime invariant divisor on $D'$. Therefore the restriction map
$$\H^0(L(-D)\vert_{D'})\longrightarrow\H^0(L(-D)\vert_x)$$
is also surjective.

Since the composition of the above maps is surjective, we get a contradiction
to the assumption that $x$ is in the base locus of $L(-D)$. 
\end{proof}

We now give the proof of the proposition for an arbitrary $l\geq 1$.

\begin{proof}[Proof of Proposition~\ref{induction}]
 We make induction on $n$, the case
$n=1$ being trivial. Note that since $l\geq 1$, $L$ is ample

Let us assume now that $n=2$. Clearly, it is enough to prove that
$(L(-D)\cdot D)\geq l-1$. Since $(L(-D)\cdot D)=(L\cdot D)-(D^2)$,
we may restrict ourselves to the case $(D^2)\geq 2$. From the
description of the selfintersection numbers in terms of the fan $\Delta$,
it follows easily that if $D'$ and $D''$ are the divisors whose rays
are adjacent to the ray corresponding to $D$, then
$(D'^2)\leq 0$ or $(D''^2)\leq 0$.

But if, for example, $(D'^2)\leq 0$, then
$L(-(l-1)D')$ is ample, so that Lemma~\ref{l1} implies that $L(-(l-1)D'-D)$
is globally generated and therefore
$$0\leq (L(-(l-1)D'-D)\cdot D)=(L(-D)\cdot D)-(l-1),$$
which completes the case $n=2$.

Suppose now that $n\geq 3$ and let $\tau\in\Delta_{n-1}$
be such that $C=V(\tau)$.
We can choose a prime invariant divisor $D'$ such that
$D'\neq D$ and $C\subset D'$.
Therefore $(L(-D)\cdot C)=(L(-D)\vert_{D'}\cdot C)$, and we may 
clearly restrict to the case when $D\cap D'\neq \emptyset$,
so that it 
is a prime invariant divisor on $D'$. We apply the induction hypothesis
for $L\vert_{D'}$; note that for every integral invariant curve
$C'\subset D'$,
$$(L\vert_{D'}\cdot C')=(L\cdot C')\geq l.$$
This concludes the proof. 
\end{proof}

We can now prove the strong form of Fujita's conjecture 
for the toric case.

\begin{proof}[Proof of Theorem~\ref{fujita}]
(i) It is clear that we may assume $n\geq 2$ and $X\not\iso\PP^n$.
We make induction on $n$.
If $L(-D_1-\ldots-D_m)$ is not globally generated, then
$$(L(-D_1-\ldots-D_m)\cdot V(\tau))<0$$
for some $\tau\in\Delta_{n-1}$. We will show that this asumption 
implies $X\iso\PP^n$, a contradiction.

We can immediately restrict ourselves to the following situation:
$2\leq m\leq n+1$, $D_1$ and $D_2$ are the divisors corresponding to the rays 
spanning together with $\tau$ maximal cones and $D_3,\ldots,D_{n+1}$ are
the divisors containing $V(\tau)$.

\smallskip

{\it Claim.\ } We have $m=n+1$, $D_i\iso\PP^{n-1}$ for every $i$,
 $1\leq i\leq n+1$, and $D_i\cap D_j\neq\emptyset$ for every $i\neq j$.

Fix $i$ such that $i\leq m$. Since $n\geq 2$,
our hypothesis and Proposition~\ref{induction} imply that $L(-D_i)$ is ample.
Hence Corollary~\ref{vanishing2} (iii) shows that the restriction map
$$\H^0(L(-D_1-\ldots-D_m))\longrightarrow
\H^0(L(-D_1-\ldots-D_m)\vert_{D_i})$$
is surjective. Since $V(\tau)\subset {\rm Bs}\,L(-D_1-\ldots-D_m)$,
it follows that $L(-D_1-\ldots-D_m)\vert_{D_i}$ is not globally 
generated.

Another application of Proposition~\ref{induction}
 gives $(L(-D_i)\cdot C)\geq n-1$
for every integral invariant curve $C\subset X$.
In particular, $(L(-D_i)\vert_{D_i}\cdot C')\geq n-1$
for every integral invariant curve $C'\subset D_i$. From the induction
hypothesis we get $D_i\iso\PP^{n-1}$, $m=n+1$ and $D_i\cap D_j\neq
\emptyset$ for $j\neq i$.

\smallskip

It is now easy to see that $X\iso\PP^n$. The claim implies that if
$D_i=V(\tau_i)$,$1\leq i\leq n+1$, and if $\tau_0$ is any other
one-dimensional cone in $\Delta$, then $\tau_0$ and $\tau_i$
do not span a cone in $\Delta$ for any $i$.
 From this it follows that the only 
one-dimensional cones in $\Delta$ are $\tau_1,\ldots,\tau_{n+1}$.
Since $X$ is smooth, it follows that $X\iso\PP^n$.

\smallskip

(ii) Since an ample line bundle on a complete smooth toric   
variety is very ample (see \cite{De}), it is enough to prove 
that if $L(-D_1-\ldots -D_m)$ is not ample, then $X\iso\PP^n$.
Again we may assume $n\geq 2$.

If $L(-D_1-\ldots-D_m)$ is not ample, then there exists
an invariant integral curve $C\subset X$ such that
$$(L(-D_1-\ldots-D_m)\cdot C)\leq 0.$$
As above, we may assume that $D_1$ and $C$ correspond
to cones in $\Delta$ spanning together a maximal cone.

By Proposition~\ref{induction}, we may apply (i) to $L(-D_1)$ and conclude that
if $X\not\iso\PP^n$, then $L(-2D_1-D_2-\ldots-D_m)$ is globally 
generated. In particular,
$$(L(-2D_1-D_2-\ldots-D_m)\cdot C)\geq 0,$$
so that
$$(L(-D_1-\ldots-D_m)\cdot C)\geq 1,$$
a contradiction.
\end{proof}

\bigskip

We conclude this section by giving two results with the same flavour 
as those proved above. By Lemma~\ref{l1}, if $L$ is ample, then $L(-D)$
is globally generated for every integral invariant divisor. The
case $X=\PP^n$, $L=\O(1)$ shows that this is optimal. The next proposition
gives the condition under which for $L$ ample we get
$L(-D_1-D_2)$ globally generated for distinct divisors $D_1$
and $D_2$ as above.

\begin{proposition}\label{next_case}
 Let $X$ be a projective smooth toric variety with
$\dim(X)=n$, $L\in\Pic(X)$ ample and $D_1$, $D_2$ distinct prime
invariant divisors. Then $L(-D_1-D_2)$ is not globally generated if
and only if 
there is an $(n-1)$-dimensional cone $\tau\in\Delta$ such that if
$\tau_1$, $\tau_2$ are the one dimensional cones correponding to
$D_1$ and $D_2$, then $(\tau, \tau_1)$ and
$(\tau, \tau_2)$ span cones in $\Delta_n$ and
$(L\cdot C)=1$, where $C=V(\tau)$.
\end{proposition}

\begin{proof}
  The ``if'' part is clear, since in this case
we have $(L(-D_1-D_2)\cdot V(\tau))=-1$.

Suppose now that $L(-D_1-D_2)$ is not globally generated. We 
prove the proposition by induction on $n$. The case $n=1$ is trivial, and
therefore we may assume $n\geq 2$. Let $x\in X$ be a fixed point 
in the base locus
of $L(-D_1-D_2)$.

Suppose first that there is an invariant prime divisor $D\neq D_1$,
$D_2$ such that $x\in D$. We apply the induction hypothesis 
for the smooth toric variety $D$, the line bundle $L\vert_D$
and the prime invariant divisors $D\cap D_1$ and $D\cap D_2$. 
By Corollary~\ref{vanishing2} (iii) , the restriction map
$$\H^0(L(-D_1-D_2))\longrightarrow\H^0(L(-D_1-D_2)\vert_D)$$
is surjective, so that our hypothesis on $x$ and $D$ implies that
$x$ is in the base locus of $L(-D_1-D_2)\vert_D$. Lemma~\ref{l1}
implies that $D\cap D_1$ and $D\cap D_2$ are nonempty. 
If $D=V(\tau_0)$, then by induction we find a cone $\tau'$ in the fan
${\rm Star(\tau_0)}$ of $D$. This corresponds to a cone $\tau\in\Delta$
which satisfies the requirements of the proposition.

Therefore it remains to consider the case when, for every fixed point $x$
in the base locus of $L(-D_1-D_2)$ and every divisor $D$ containing
$x$, we have $D=D_1$ or $D=D_2$. Clearly this implies $n=2$
and the fact that the base locus consists of a point, the corresponding
cone being generated by the rays defining $D_1$ and $D_2$. But this 
contradicts Corollary~\ref{base_locus} and the proof is complete.
\end{proof}

As a consequence of Proposition~\ref{induction}
 we get that if $(L\cdot C)\geq 2$
for every integral invariant curve on $X$, then $L(-D)$ is ample for
every prime invariant divisor $D$. The next result makes this more precise
by giving the condition for an ample line bundle $L$ and a prime invariant
divisor $D$ to have $L(-D)$ not ample.

\begin{proposition}\label{ample_case}
 Let $X$ be a complete smooth toric variety
with $\dim(X)=n$, $L\in\Pic(X)$ an ample line bundle and 
$D=V(\tau_0)$ a prime invariant divisor. Then $L(-D)$ is not ample
if and only if
 there is $\tau\in\Delta_{n-1}$ such that $\langle\tau,\tau_0\rangle
\in\Delta
_n$ and $(L\cdot V(\tau))=1$.
\end{proposition}

\begin{proof}
It is clear that if there exists $\tau$
as above, then $(L(-D)\cdot V(\tau))=0$, so that $L(-D)$
is not ample.

Suppose now that $L(-D)$ is not ample and therefore there exists
$\tau'\in\Delta_{n-1}$ such that $(L(-D)\cdot V(\tau'))\leq 0$.
Since $L(-D)$ is globally generated by Lemma~\ref{l1}, we must have
$(L(-D)\cdot V(\tau'))=0$.

We must have $(D\cdot V(\tau'))\neq 0$,
and therefore we deduce that either
$\langle\tau_0,\tau'\rangle\in\Delta_n$ or $V(\tau')\subset D$.
In the first case, we have $(L\cdot V(\tau'))=1$ and may take
$\tau=\tau'$.

If $V(\tau')\subset D$, we choose a divisor $D_1=V(\tau_1)$
such that $\langle\tau_1, \tau'\rangle\in\Delta_n$. Then
$(L(-D-D_1)\cdot V(\tau')) < 0$, and Proposition~\ref{next_case}
 implies that there
is $\tau\in\Delta_{n-1}$ such that
 $\langle\tau_0, \tau\rangle\in\Delta_n$
and $(L\cdot V(\tau))=1$.
\end{proof}

\section{Sections of ample line bundles}

In this section we fix a globally generated line bundle $L$
on a complete toric variety $X$ and an invariant
divisor $D$ such that $L\iso\O(D)$. Since $L$ is globally generated,
for each
maximal cone $\sigma$ there is a unique $u_{\sigma}\in M$
such that ${\rm div}(\chi^{u_{\sigma}})+D$ is effective and 
zero on $U_{\sigma}$. Equivalently, for each maximal cone $\sigma$,
there is a nonzero section $s_{\sigma}\in\H^0(X, L)$,
unique up to scalars, which is an eigenvector with respect to the
torus action and whose restriction to $U_{\sigma}$ is everywhere 
nonzero.

A well-known ampleness criterion (see \cite{Fu}, Section 3.4) 
says that $L$ is ample if and only if $u_{\sigma}\neq u_{\tau}$
(or, equivalently, $k\,s_{\sigma}\neq k\,s_{\tau}$) for 
$\sigma\neq\tau$. From the unicity of the sections $s_{\sigma}$,
this is equivalent to the fact that if $\sigma\neq\tau$, then
$s_{\sigma}\vert_{U_{\tau}}$ vanishes at some point. But in that case,
it must vanish at the unique fixed point $x_{\tau}$
of  $U_{\tau}$.

We consider the following map whose components are given by the
restriction maps:
$$\phi : \H^0(L)\longrightarrow\bigoplus_{\sigma\in\Delta_{\rm max}}
\H^0(L\vert_{\{x_{\sigma}\}}).$$
Since $\phi$ is an equivariant map under the torus action, the discussion
above shows that $L$ is ample if and only if $\phi$ is surjective.

\smallskip

Our goal in this section is to extend this property
of ample line bundles in the case when $X$
is smooth to a set of higher dimensional subvarieties which are
pairwise disjoint. More precisely, we have the following

\begin{theorem}\label{surjectivity}
 Let $X$ be a projective smooth toric variety
and $L\in\Pic(X)$ an ample line bundle. If $Y_1,\ldots,Y_r\subset X$
are integral invariant subvarieties such that $Y_i\cap Y_j=\emptyset$
for $i\neq j$ and
$$\psi : \H^0(L)\longrightarrow\bigoplus_{i=1}^r
\H^0(L\vert_{Y_i})$$
is induced by restrictions, then $\psi$ is surjective.
\end{theorem}

\begin{proof}
Let $Y=\bigcup_{i=1}^r Y_i$. In order to prove
that
$$\psi : \H^0(L)\longrightarrow \H^0(L\vert_Y)$$
is surjective, it is enough to prove that $\H^1(L\tensor {\mathcal I}
_{Y/X})=0$.

Let $\pi : \widetilde{X}\longrightarrow X$ be the blowing-up of $X$
along $Y$ and $E$ the exceptional divisor. Then
$$\H^1(X, L\tensor {\mathcal I}_{Y/X})
\iso\H^1(\widetilde{X}, \pi^*L\tensor\O(-E)).$$
Since $X$ is smooth, the blowing-up of $X$ along an integral
invariant subvariety is still a smooth toric variety (\cite{Ew}).
Since $Y_i\cap Y_j=\emptyset$ for $i\neq j$, $\pi$ is a composition
of such transformations, and therefore $\widetilde{X}$ is a toric 
variety. Moreover, from the description in Ewald [1996] it
follows that if $E_i=\pi^{-1}(Y_i)$, then $E_i$ is an invariant
prime divisor on $\widetilde{X}$ and $E=\sum_{i=1}^r E_i$.

Since $L$ is ample, Proposition 7.10 in \cite{Ha} implies
that there is an integer $s\geq 1$ such that $\pi^*(L^s)\tensor\O(-E)$
is ample on $\widetilde{X}$.
We choose an invariant divisor $D$ on $\widetilde{X}$ such that
$\pi^*L\iso\O(D)$. Then $D-(1/s)E=D-(1/s)\sum_{i=1}^r E_i$
is $\QQ$-ample and $\lfloor D-(1/s)E\rfloor=D-E$.
Now Corollary~\ref{vanishing2} gives $\H^1(X, \pi^*L\tensor\O(-E))=0$.
\end{proof}

%----------------------------------------
%%  REFERENCES
%%----------------------------------------------------------
%%\bibliographystyle{amsalpha}
%%\bibliography{MRC}

\begin{thebibliography}{BC}



\bibitem[BC]{BC}
V. Batyrev and D. Cox, \emph{ On the Hodge structure of projective
  hypersurfaces in toric varieties}, Duke Math. J. 75 
(1994), 293--338. 

\bibitem[BTLM]{BTLM}
A. Buch, J. F. Thomsen, N. Lauritzen, V. Mehta,
\emph{The Frobenius morphism of a toric variety}, Toh\^{o}ku Math. J. 
49 (1997), 355--366.  

\bibitem[Cox]{Co}
D. Cox, \emph{The homogeneous coordinate ring of a toric variety},
 J. Algebraic Geom. 4 (1995), 17--50.

\bibitem[De]{De}
M. Demazure, \emph{Sous-groupes alg\'ebriques de rang maximum
   du groupe de Cremona}, Ann. Sci. Ecole Norm. Sup. 3
   (1970), 507--588.

\bibitem[DR]{DR}
S. Di Rocco, \emph{Generation of $k$-jets on toric varieties},
    Math. Z. 231(1999), 169--188.

\bibitem[EL]{EL}
L. Ein and R. Lazarsfeld, \emph{Syzygies and Koszul cohomology
   of smooth projective varieties of arbitrary dimension},
   Invent. Math. 111 (1993), 51--67.

\bibitem[Ei]{Ei}
D. Eisenbud, \emph{Commutative algebra with a view toward
   algebraic geometry}, Grad. Texts in Math. 150, Springer,
   New York, 1995.

\bibitem[EMS]{EMS}
D. Eisenbud, M. Musta\c{t}\v{a} and M. Stillman, \emph{Cohomology
   on toric varieties and local cohomology with
   monomial supports}, J. Symbolic Comput. 29 (2000),
   583--600.

\bibitem[Ew]{Ew}
G. Ewald, \emph{Combinatorial convexity and algebraic geometry},
   Grad. Texts in Math. 168, Springer, New York, 1996.

\bibitem[Fu]{Fu}
W. Fulton, \emph{Introduction to toric varieties}, Ann. of Math. Stud.
    131, The William H.~Rover Lectures in Geometry,
   Princeton Univ. Press, Princeton, NJ, 1993.

\bibitem[Ha]{Ha}
R. Hartshorne, \emph{Algebraic geometry}, 
   Grad. Texts in Math. 52, Springer-Verlag, New York, 1977.

\bibitem[Ka]{Ka}
T. Kajiwara, \emph{The functor of a toric variety with enough
   effective invariant divisors}, Toh\^{o}ku Math. J. 50 (1998),
   139--157.

\bibitem[KKMS]{KKMS}
G. Kempf, F. Knudsen, D. Mumford, B. Saint-Donat, \emph{Toroidal
embeddings I}, Lecture Notes in Math. 339, Springer-Verlag,
Berlin-New York, 1973.

\bibitem[La]{La}
R. Lazarsfeld, \emph{Lectures on linear series}, with the assistance
   of Fern\'andez del Busto, IAS/Park City Math. Ser. 3, Complex
   algebraic geometry (Park City, UT, 1993),
    161--219, Amer. Math. Soc. Providence, RI, 1997.

\bibitem[Ma]{Ma}
A. Mavlyutov, \emph{Semi-ample hypersurfaces in toric varieties},
   Duke Math. J. 101 (2000), 85--116.

\bibitem[MR]{MR}
V. B. Mehta and A. Ramanathan, \emph{Frobenius splitting
   and cohomology vanishing for Schubert varieties}, Ann. of Math.(2)
   122 (1985), 27--40.

\bibitem[Mu]{Mu}
M. Musta\c{t}\v{a}, \emph{Local cohomology at monomial ideals},
   J. Symbolic Comput. 29 (2000), 709--720.

\bibitem[Oda]{Od}
T. Oda, \emph{Convex bodies and algebraic geometry},
   Ergeb. Math. Grenzgeb. (3) 15, 
   Springer-Verlag, Berlin Heidelberg New York, 1988.

\bibitem[Re]{Re}
M. Reid, \emph{Decomposition of toric morphisms}, in
   Arithmetic and Geometry II, Progr. Math. 36, Birkh\"{a}user
   Boston, MA, 395--418.

\bibitem[Sm]{Sm}
K. Smith, \emph{Fujita's freeness conjecture in terms of
   local cohomology}, J.Algebraic Geom. 6 (1997), 417--429.

\bibitem[Ya]{Ya}
K. Yanagawa, \emph{Alexander duality for Stanley-Reisner rings and
   square-free ${\NN}^n$-graded modules}, J. Algebra 225 (2000),
   630--645.

\end{thebibliography}

\providecommand{\bysame}{\leavevmode\hbox to3em{\hrulefill}\thinspace}

\end{document}